\documentclass[11pt,reqno]{amsart}

\usepackage{amsmath,amsthm,amssymb,enumitem}

\newcommand{\Z}{\mathbb{Z}}

\newcommand{\R}{\mathbb{R}}

\newcommand{\innerProd}[2]{\langle#1,#2\rangle}

\newcommand{\trace}{\ensuremath{\text{Tr}}}
\newcommand{\diag}{\ensuremath{\text{diag}}}

\newtheorem{proposition}{Proposition}
\newtheorem{theorem}{Theorem}

\theoremstyle{definition}

\theoremstyle{remark}

\usepackage{algorithm,algorithmicx}
\usepackage{algpseudocode}
\usepackage{bbm}
\usepackage[mathlines,pagewise,left]{lineno}
\usepackage{amssymb}
\usepackage{color,xcolor,tikz}
\usetikzlibrary{arrows,shapes,chains}  
\usepackage{amsmath}
\usepackage{./bcol}
\usepackage{multirow, array}
\usepackage{caption}
\usepackage{xspace}
\usepackage[numbers]{natbib}
\usepackage{comment}
\usepackage{graphicx} 
\usepackage{booktabs} 
\usepackage{pdflscape}
\usepackage[export]{adjustbox}
\usepackage[font=scriptsize]{subfig}
\usepackage{mathtools}
\usepackage{hyperref}
\usepackage[titletoc]{appendix}

\usepackage{amsthm}
\theoremstyle{definition}
\theoremstyle{remark}

\newcommand{\be}[1]{\begin{equation}\label{#1}}
\newcommand{\ee}{\end{equation}}

\newcommand{\sdpstandard}{\ensuremath{{(\textsf{Shor})}}\xspace}

\newcommand{\optPersp}{\ensuremath{{(\textsf{OptPersp})}}\xspace}

\newcommand{\CQI}{\ensuremath{{(\textsf{QI})}}\xspace}

\newcommand{\set}[1]{\ensuremath{\mathcal{#1}}}

\newcolumntype{\resetRow}{>{\global\let\currentrowstyle\relax}}
\newcolumntype{^}{>{\currentrowstyle}}

\def\SingleSpacedXI{\linespread{1.1}}
\SingleSpacedXI

\usepackage[colorinlistoftodos,prependcaption,textsize=small]{todonotes}
\title[Perspective reformulation v.s. Shor's SDP for MIQP]{The equivalence of optimal perspective formulation and Shor's SDP for quadratic programs with indicator variables}
\author{Shaoning Han, Andr\'{e}s G\'{o}mez and Alper Atamt\"urk}

\thanks{ \noindent \hskip -5mm
	S. Han, A. G\'{o}mez: Daniel J. Epstein Department of Industrial and Systems Engineering, Viterbi School of Engineering, University of Southern California, CA 90089. \texttt{gomezand@usc.edu}, \texttt{shaoning@usc.edu}\\
A. Atamt\"urk: Department of Industrial Engineering \& Operations Research, University of California, Berkeley, CA 94720.
\texttt{atamturk@berkeley.edu}  
}
\begin{document}
	\maketitle
	
	\begin{abstract}
		\vskip 3mm
		\noindent 
In this paper, we compare the strength of the optimal perspective reformulation and Shor's SDP relaxation. We prove that these two formulations are equivalent for quadratic optimization problems with indicator variables. \\
		 
\end{abstract}

\begin{center}
	{December} 2021 \\
\end{center}

\section{Introduction}

\label{sec:intro}

We consider the optimization problem with a convex quadratic objective with indicators:
\begin{subequations}\label{eq:mixed_integer_problem}
	\begin{align*}
		\min\;& a'x+b'y+ y'Qy \\ 
		\CQI \qquad\text{s.t.}\;&y_i(1-x_i)=0&\forall i\in[n]\\
		&x\in\{0,1\}^n,\;y\in\R_+^n\\
		&{(x,y)\in \set{X}\subseteq \R^n\times\R^n,}
	\end{align*}
\end{subequations}
where $a$ and $b$ are $n$-dimensional vectors, $Q\in\R^{n\times n}$ is a {symmetric} positive semidefinite (PSD) matrix, 
$[n] := \{1,2, \ldots, n\}$,  {and \set X is a closed convex set that incorporates other unspecified constraints arising in applications.} For each $i \in [n]$, the {complementarity} constraint $y_i (1-x_i) = 0$, along with the indicator variable $x_i \in \{0,1\}$, is used to state that $y_i=0$ whenever $x_i=0$. 
Numerous applications, including portfolio optimization \cite{bienstock1996computational}, optimal control  \cite{gao2011cardinality}, image segmentation  \cite{hochbaum2001efficient}, signal denoising \cite{bach2019submodular}
are either formulated as \CQI or can be relaxed to \CQI. While special cases of \CQI with a diagonal $Q$ or M-matrix $Q$ can be solved in polynomial time \cite{atamturk2018strong}, in general, \CQI is \NP-hard \cite{HGA:2x2}.

Building strong convex relaxations of \CQI 
is instrumental in solving it effectively. There is a substantial body of research on the perspective formulation of convex univariate functions with indicators \cite{akturk2009strong,dong2015regularization,dong2013,frangioni2006perspective, gunluk2010perspective,
	hijazi2012mixed,wu2017quadratic}. When $Q$ is diagonal, $y'Qy$ is separable and the perspective formulation provides the convex hull of the epigraph of $y'Qy$ with indicator variables by strengthening each term $Q_{ii}y_i^2$ with its perspective counterpart $Q_{ii}y_i^2/x_i$, individually.

 Another powerful approach for nonconvex quadratic optimization problems is semidefinite programming (SDP) reformulation,
 first proposed by Shor \cite{shor1987quadratic}. 
 Specifically, a convex relaxation is constructed by introducing a rank-one matrix $Z$ representing $zz'$, where $z$ is the decision vector, and then forming the semidefinite relaxation $Z \succeq zz'$. Such SDP relaxations have been widely utilized in numerous applications, including max-cut problems \cite{goemans1995improved}, hidden partition problems of finding clusters in large network datasets \cite{javanmard2016phase},  matrix completion problems \cite{alfakih1999solving,candes2010matrix}, power systems \cite{FALA:conic-uc}, 
 robust optimization problems \cite{ben2009robust}.
  Sufficient conditions for exactness of SDP relaxations \cite[e.g.][]{burer2019exact,ho2017second,jeyakumar2014trust, wang2019tightness,wang2019generalized} {and stronger rank-one conic formulations \cite{AG:rank-one,AG:supermodular}} are also {given} in the literature.

These two approaches have been studied extensively in literature. While it is known that Shor's SDP formulation is at least as strong as the optimal perspective formulation \cite{dong2015regularization}, the other direction has not been explored. We show in this note that these two formulations are, in fact, equivalent. The equivalence makes the perspective formulation the favorable choice as it is much smaller than Shor's SDP and easier to solve. 

\subsection*{Notation} Throughout, we adopt the following convention for division by $0$: given $a\geq 0$, $a^2/0=\infty$ if $a\neq 0$ and $a^2/0=0$ if $a=0$.
For a vector $v$, $\diag(v)$ denotes the diagonal matrix $V$ with $V_{ii} = v_i$ for each $i$. 
\section{Optimal perspective formulation vs. Shor's SDP}\label{sec:previousResults}
In this section we analyze two well-known convex formulations: the optimal perspective formulation and Shor's SDP. We first introduce the two formulations and then show that they are equivalent for \CQI.
 
Splitting $Q$ into some diagonal PSD matrix $D=\diag(d)$ and a PSD residual, i.e., $Q-D\succeq 0$, one can apply the perspective reformulation to each diagonal term, by replacing $D_{ii}y_i^2$ with $D_{ii}y_i^2/x_i$,
to get a valid convex relaxation of \CQI --after relaxing integrality constraints in $x$ and dropping the {complementarity} constraints $y_i(1-x_i)=0$:
\begin{equation}\label{eq:normal_persp}
	\begin{aligned}
		\min_{(x,y)\in\R^{2n}}\;&a'x+b'y+y'(Q-\diag(d))y+\sum_{i\in[n]}d_i\frac{y_i^2}{x_i}\\
		\text{s.t. }&0\le x\le1,\,y\ge0\\
		&{(x,y)\in \set{X}\subseteq \R^n\times\R^n.}
	\end{aligned}
\end{equation}
{In certain applications, such as the sensor placement \cite{frangioni2011projected}, the single-period unit commitment \cite{galiana2003reconciling} and the $\ell_2$-penalized least square regression \cite{pilanci2015sparse}, vector $d$ is immediate from the context. Thus, in such cases, \eqref{eq:normal_persp} directly delivers a strong relaxation of \eqref{eq:mixed_integer_problem}. For cases where a decomposition of $Q$ is not immediate, several approaches are proposed in literature \citep[e.g.][]{frangioni2007sdp,zheng2014improving} to obtain a desirable $d$. Because different decompositions  usually yield different relaxations, the relaxation quality of \eqref{eq:normal_persp} relies on the choice of vector $d$.} Introducing a symmetric matrix variable $Y$,
\citet{dong2015regularization} describe an optimal perspective relaxation for \CQI: 
\begin{subequations}\label{eq:perspective}
		\begin{align}
			\min\;& a'x+b'y+\innerProd{Q}{Y}\\
			\optPersp\qquad\qquad\text{s.t.}\;&Y-yy'\succeq 0\label{eq:perspective_psd}\\
			&y_i^2\leq Y_{ii}x_i&\forall i\in[n]\label{eq:perspective_rotated}\\
			&0\le x \le 1,\;y\ge0\\
			&{(x,y)\in \set{X}\subseteq \R^n\times\R^n.}
		\end{align}
\end{subequations}
{Note that by adding integrality constraints on $x$ in \optPersp, one obtains a mixed-integer SDP problem, which can be solved by a branch-and-bound algorithm. What is more, the resulting mixed-integer program is equivalent to the original model \CQI. Indeed, if $x$ is integral, then  \eqref{eq:perspective_rotated} either reduces to $y_i=0$ or is implied \eqref{eq:perspective_psd}. In any case, $Y$ and $y$ are linked only through \eqref{eq:perspective_psd}, which must hold as an equality at the optimal solution because $Q\succeq 0$.} 

The authors \cite{dong2015regularization} show that \optPersp is optimal in the sense that every perspective relaxation of the form \eqref{eq:normal_persp} is dominated by \optPersp, and moreover, there indeed exists a decomposition of $Q$ such that the resulting perspective formulation \eqref{eq:normal_persp} is equivalent to \optPersp.
\begin{proposition}[Theorem~3 in \cite{dong2015regularization}]\label{prop:decomp-diag}
	\optPersp is equivalent to the {following} max-min optimization {problem}:
	\begin{subequations}
		\begin{align*}
			\max_{d\in\R^n}\min_{(x,y)\in\R^{2n}}\;&a'x+b'y+y'(Q-\diag(d))y+\sum_{i\in[n]}d_i\frac{y_i^2}{x_i}\\
			\text{s.t. }
			&{d \ge 0,}	\;Q-\diag(d)\succeq0\\
			&0\le x\le1,\,y\ge0\\
			&{(x,y)\in \set{X}\subseteq \R^n\times\R^n.}
		\end{align*}
	\end{subequations}
\end{proposition}
\ignore{
}

Next, we consider Shor's SDP relaxation for problem~\CQI: 
	\begin{subequations}
		\begin{align}
			\min\;& a'x+b'y+\sum_{i=1}^n\sum_{j=1}^n Q_{ij}Z_{ij}\\
			\text{s.t.}\;&y_i-Z_{i,i+n}=0 & \forall i\in[n]\\
			\sdpstandard\qquad \qquad &x_i-Z_{i+n,i+n}=0&\forall i\in[n]\\
			&Z-\begin{pmatrix}
				y\\x
			\end{pmatrix}\begin{pmatrix}
				y'&x'
			\end{pmatrix}\succeq 0\label{eq:sdp_psd}\\
			&0\le x \le 1,\;y\ge0\\
			&{(x,y)\in \set{X}\subseteq \R^n\times\R^n,}
		\end{align}
	\end{subequations}
	where $Z\in \R^{2n\times 2n}$ such that  $Z_{ii}$ is a proxy for $y_i^2$,  $Z_{i+n,i+n}$ is a proxy for $x_i^2$, $Z_{i,i+n}$ is a proxy for $x_iy_i$, $i\in[n]$, and $Z_{ij}$ is a proxy for $y_iy_j$ for $1\leq i,j\leq n$. 
	
	It is known that \sdpstandard\ is at least as strong as \optPersp \cite{dong2013}, as constraints \eqref{eq:perspective_rotated} are implied {by} the positive definiteness of some $2\times 2$ principal minors of \eqref{eq:sdp_psd}. We show below that the two formulations are, in fact, equivalent. As \optPersp is a much smaller formulation than \sdpstandard, the equivalence makes it the favorable choice. 
	\begin{theorem}\label{thm:equivalence}
	\optPersp is equivalent to 	\sdpstandard.
	\end{theorem}
	\begin{proof}
		First we verify \sdpstandard\ is at least as strong as \optPersp by checking that constraints \eqref{eq:perspective_psd}--\eqref{eq:perspective_rotated} are implied {by} \eqref{eq:sdp_psd}.  
		Let $Y_{ij} = Z_{ij}$ for any ${i,j\in[n]}$. By Schur Complement Lemma,
		\[ Z-\begin{pmatrix}
			y\\x
		\end{pmatrix}\begin{pmatrix}
			y'&x'
		\end{pmatrix}\succeq 0\iff \begin{pmatrix}
			1&\begin{matrix}
				y'&x'
			\end{matrix}\\
			\begin{matrix}
				y\\x
			\end{matrix}&Z
		\end{pmatrix}\succeq 0. \]
		Since $Y$ is a principle submatrix of $Z$, we have \[ \begin{pmatrix}
			1&y'\\y&Y
		\end{pmatrix} \succeq0 \Leftrightarrow Y-yy'\succeq 0.\] Moreover, constraint \eqref{eq:sdp_psd} also implies that for any ${i\in[n]}$,
		\begin{equation}\label{eq:psdSubmatrices}
			\begin{pmatrix}
				Z_{ii}&Z_{i,i+n}\\Z_{i,i+n}&Z_{i+n,i+n}
			\end{pmatrix}\succeq 0.
		\end{equation} 
		After substituting $Y_{ii} = Z_{ii}, \ x_i =  Z_{i+n,i+n},$ and $y_i = Z_{i,i+n}$, we find that \eqref{eq:psdSubmatrices} {implies} $Y_{ii}x_i\ge y_i^2$ in \optPersp, concluding the argument. 
		
		We next prove that \optPersp is at least as strong as \sdpstandard, by showing that for any given feasible solution $(x,y,Y)$ of  \optPersp, it is possible to construct a feasible solution of \sdpstandard\ with same objective value. First, we rewrite $Z$ in the form  $Z=\begin{pmatrix}
			Y&U\\U'&V
		\end{pmatrix}.$ 
		For a fixed feasible solution $(x,y,Y)$ of \optPersp consider the optimization problem
		\begin{subequations}\label{eq:auxPrimal}
			\begin{align}
				\lambda^*:=&\min_{\lambda,U,V}\;\lambda\\
				\text{s.t. }&\begin{pmatrix}
					1&y'&x'\\
					y&Y &U\\
					x&U'&V
				\end{pmatrix}+\lambda I\succeq 0 \label{eq:primal_matrix_inequality}\\
				&U_{ii}=y_i,&\forall i\in[n] \,\, \\
				&V_{ii}=x_i,&\forall i\in[n],
			\end{align}
		\end{subequations}
		where $I$ is the  identity matrix. Observe that if $\lambda^*\le0$, then an optimal solution of \eqref{eq:auxPrimal} satisfies \eqref{eq:sdp_psd} and thus induces a feasible solution of \sdpstandard. We show next that this is, indeed, the case.
		
        Let $\tilde{Y}=\begin{pmatrix}
			1&y'\\y&Y
		\end{pmatrix}$ and consider the SDP dual of \eqref{eq:auxPrimal}:
		\begin{subequations}
			\begin{align} \label{eq:auxDual}
				\lambda^*=\max_{R,s,t,z} \;&-\innerProd{\tilde{Y}}{R} - \sum_i( 2x_iz_i+t_ix_i+2s_iy_i) &\\
				\text{s.t. }& \begin{pmatrix}
					R&\begin{matrix}
						z'\\\diag(s)
					\end{matrix}\\\label{eq:dual_variable}
					\begin{matrix}
						z&\diag(s)
					\end{matrix}&\diag(t)
				\end{pmatrix}\succeq 0\\
				&\trace(R) + \sum_i t_i =1,&\label{eq:normalization}
			\end{align}
		\end{subequations}
		where $R, z, \diag(s), \diag(t)$ are the dual variable associated with $\tilde{Y}+\lambda I, x, U,$ and $V+\lambda I$, respectively. Note that we abuse the symbol $I$ to represent the identity matrices of different dimensions. 
				One can verify that the strong duality holds for \eqref{eq:auxPrimal} since $\lambda$ can be an arbitrary positive number to ensure that the matrix inequality holds strictly.
		Because the off-diagonal elements of $U$ and $V$ do not appear in the primal objective function and constraints  other than \eqref{eq:primal_matrix_inequality}, the corresponding dual variables are zero.
		
		Note that to show $\lambda^*\le0$,  it is sufficient to consider a relaxation of \eqref{eq:auxDual}. Therefore, dropping
		\eqref{eq:normalization}, 
		it is sufficient to show that
		\begin{equation}
			\innerProd{\tilde{Y}}{R} + \sum_i( 2x_iz_i+t_ix_i+2s_iy_i)\geq 0\label{eq:objective}
		\end{equation} for all $t\ge0$, $s$, $z,$ and $R$ satisfying \eqref{eq:dual_variable}. 
		
		Observe that if $t_i=0$, then $s_i=z_i=0$ in any solution satisfying \eqref{eq:dual_variable}. In this case, all such terms indexed by $i$ vanish in \eqref{eq:objective}. Therefore, it suffices to prove $\eqref{eq:objective}$ holds for all $t>0$.
		
		For $t>0$, by Schur Complement Lemma, \eqref{eq:dual_variable} is equivalent to
		\begin{equation}
			R\succeq \begin{pmatrix}
				z'\\\diag(s)
			\end{pmatrix}\diag^{-1}(t)\begin{pmatrix}
				z&\diag(s)
			\end{pmatrix}.\label{eq:dual_variable_equiv}
		\end{equation}
		Moreover, since $\tilde{Y}\succeq0$, we find that $$\innerProd{\tilde{Y}}{R}\ge\innerProd{\tilde{Y}}{\begin{pmatrix}
				z'\\\diag(s)
			\end{pmatrix}\diag^{-1}(t)\begin{pmatrix}
				z&\diag(s)
		\end{pmatrix}}$$
		whenever $R$ satisfies \eqref{eq:dual_variable_equiv}. Substituting the term $\innerProd{\tilde{Y}}{R}$ in \eqref{eq:objective} by its lower bound, it suffices to show that 	
		\begin{equation}\label{eq:lb_objecive}
			\innerProd{\tilde{Y}}{\begin{pmatrix}
					z'\\\diag(s)
				\end{pmatrix}\diag^{-1}(t)\begin{pmatrix}
					z&\diag(s)
			\end{pmatrix}} + \sum_i( 2x_iz_i+t_ix_i+2s_iy_i)\ge0,
		\end{equation} 
		holds for all $t>0,s, z$ and $R$ satisfying \eqref{eq:dual_variable_equiv}. A direct computation shows that
		\[\begin{pmatrix}
			z'\\\diag(s)
		\end{pmatrix}\diag^{-1}(t)\begin{pmatrix}
			z&\diag(s)
		\end{pmatrix}=\begin{pmatrix}
			\sum_i z_i^2/t_i&z_1s_1/t_1&\cdots&z_ns_n/t_n\\
			z_1s_1/t_1&s_1^2/t_1&&\\
			\vdots&&\ddots&\\
			z_ns_n/t_n&&&s_n^2/t_n
		\end{pmatrix}\]
		with all off-diagonal elements equal to $0$, except for the first row/column.
		Thus, \eqref{eq:lb_objecive} reduces to the separable expression
		\[ \sum_i \left( \frac{z_i^2}{t_i}+\frac{2z_is_iy_i}{t_i}+\frac{s_i^2}{t_i}Y_{ii}+2x_iz_i+x_it_i+2y_is_i \right)\ge0.  \]
		For each term, we have
		\begin{align*}
			&\frac{z_i^2}{t_i}+\frac{2z_is_iy_i}{t_i}+\frac{s_i^2}{t_i}Y_{ii}+2x_iz_i+x_it_i+2y_is_i\\
			= \ &(z_i^2+2z_is_iy_i+s_i^2Y_{ii})/t_i+x_it_i+2x_iz_i+2y_is_i\\
			\ge \ & 2\sqrt{x_i(z_i^2+2z_is_iy_i+s_i^2Y_{ii})}+2x_iz_i+2y_is_i\\
			\ge \ &0,
		\end{align*}
		where the first inequality follows from {the inequality between the arithmetic and geometric mean} $a+b\ge 2\sqrt{ab}$ for $a,b\ge0$.
		The last inequality holds trivially if $2x_iz_i+2y_is_i\geq 0$;  otherwise, we have
		\begin{align*}
			&\sqrt{x_i(z_i^2+2z_is_iy_i+s_i^2Y_{ii})}\ge -(x_iz_i+y_is_i)\\
			\iff&x_i(z_i^2+2z_is_iy_i+s_i^2Y_{ii})\ge(x_iz_i+y_is_i)^2\\
			\iff&x_iz_i^2(1-x_i)+s_i^2(x_iY_{ii}-y_i^2)\ge 0. \tag{as $0\le x_i\le 1$ and $x_iY_{ii}\ge y_i^2$}
		\end{align*}
		In conclusion, $\lambda^*\le0$ and this completes the proof.
	\end{proof}

{
	\section{Conclusion}
	In this short note, we prove that Shor's SDP and optimal perspective reformulation have the equal strength for convex quadratic optimization problems with indicator variables. The smaller size of the optimal perspective make it advantages to solve as a relaxation compared to Shor's SDP. For solving \CQI with integral constraints, it may be preferable to use conic quadratic formulation \eqref{eq:normal_persp} with an optimal $d^*$, which can be retrieved as optimal dual variables of \eqref{eq:perspective} with respect to constraints \eqref{eq:perspective_rotated}.
}

\section*{Acknowledgments}
Andr\'es G\'omez is supported, in part, by grants 1930582 and 1818700 from the National Science Foundation. Alper Atamt\"urk is supported, in part, {by} NSF AI Institute for Advances in Optimization Award 2112533, NSF grant 1807260, and DOD ONR grant 12951270.

\linespread{1.0}

\bibliographystyle{apalike}
\bibliography{./Bibliography}

\begin{thebibliography}{}

\bibitem[Akt{\"u}rk et~al., 2009]{akturk2009strong}
Akt{\"u}rk, M.~S., Atamt{\"u}rk, A., and G{\"u}rel, S. (2009).
\newblock A strong conic quadratic reformulation for machine-job assignment
  with controllable processing times.
\newblock {\em Operations Research Letters}, 37(3):187--191.

\bibitem[Alfakih et~al., 1999]{alfakih1999solving}
Alfakih, A.~Y., Khandani, A., and Wolkowicz, H. (1999).
\newblock Solving euclidean distance matrix completion problems via
  semidefinite programming.
\newblock {\em Computational Optimization and Applications}, 12(1-3):13--30.

\bibitem[Atamt{\"u}rk and G{\'o}mez, 2018]{atamturk2018strong}
Atamt{\"u}rk, A. and G{\'o}mez, A. (2018).
\newblock Strong formulations for quadratic optimization with {M}-matrices and
  indicator variables.
\newblock {\em Mathematical Programming}, 170(1):141--176.

\bibitem[Atamt\"urk and G\'omez, 2019]{AG:rank-one}
Atamt\"urk, A. and G\'omez, A. (2019).
\newblock Rank-one convexification for sparse regression.
\newblock {\em arXiv preprint arXiv:1901.10334}.

\bibitem[Atamt\"urk and G\'omez, 2020]{AG:supermodular}
Atamt\"urk, A. and G\'omez, A. (2020).
\newblock Supermodularity and valid inequalities for quadratic optimization
  with indicators.
\newblock {\em arXiv preprint arXiv:2012.14633}.

\bibitem[Bach, 2019]{bach2019submodular}
Bach, F. (2019).
\newblock Submodular functions: from discrete to continuous domains.
\newblock {\em Mathematical Programming}, 175(1-2):419--459.

\bibitem[Ben-Tal et~al., 2009]{ben2009robust}
Ben-Tal, A., El~Ghaoui, L., and Nemirovski, A. (2009).
\newblock {\em Robust Optimization}, volume~28.
\newblock Princeton University Press.

\bibitem[Bienstock, 1996]{bienstock1996computational}
Bienstock, D. (1996).
\newblock Computational study of a family of mixed-integer quadratic
  programming problems.
\newblock {\em Mathematical Programming}, 74(2):121--140.

\bibitem[Burer and Ye, 2019]{burer2019exact}
Burer, S. and Ye, Y. (2019).
\newblock Exact semidefinite formulations for a class of (random and
  non-random) nonconvex quadratic programs.
\newblock {\em Mathematical Programming}, pages 1--17.

\bibitem[{Candes} and {Plan}, 2010]{candes2010matrix}
{Candes}, E.~J. and {Plan}, Y. (2010).
\newblock Matrix completion with noise.
\newblock {\em Proceedings of the IEEE}, 98(6):925--936.

\bibitem[Dong et~al., 2015]{dong2015regularization}
Dong, H., Chen, K., and Linderoth, J. (2015).
\newblock Regularization vs. relaxation: A conic optimization perspective of
  statistical variable selection.
\newblock {\em arXiv preprint arXiv:1510.06083}.

\bibitem[Dong and Linderoth, 2013]{dong2013}
Dong, H. and Linderoth, J. (2013).
\newblock On valid inequalities for quadratic programming with continuous
  variables and binary indicators.
\newblock In Goemans, M. and Correa, J., editors, {\em Proceedings of IPCO
  2013}, page 169–180, Berlin. Springer.

\bibitem[Fattahi et~al., 2017]{FALA:conic-uc}
Fattahi, S., Ashraphijuo, M., Lavaei, J., and Atamt{\"u}rk, A. (2017).
\newblock Conic relaxations of the unit commitment problem.
\newblock {\em Energy}, 134:1079--1095.

\bibitem[Frangioni and Gentile, 2006]{frangioni2006perspective}
Frangioni, A. and Gentile, C. (2006).
\newblock Perspective cuts for a class of convex 0--1 mixed integer programs.
\newblock {\em Mathematical Programming}, 106(2):225--236.

\bibitem[Frangioni and Gentile, 2007]{frangioni2007sdp}
Frangioni, A. and Gentile, C. (2007).
\newblock {SDP} diagonalizations and perspective cuts for a class of
  nonseparable miqp.
\newblock {\em Operations Research Letters}, 35:181--185.

\bibitem[Frangioni et~al., 2011]{frangioni2011projected}
Frangioni, A., Gentile, C., Grande, E., and Pacifici, A. (2011).
\newblock Projected perspective reformulations with applications in design
  problems.
\newblock {\em Operations research}, 59(5):1225--1232.

\bibitem[Galiana et~al., 2003]{galiana2003reconciling}
Galiana, F.~D., Motto, A.~L., and Bouffard, F. (2003).
\newblock Reconciling social welfare, agent profits, and consumer payments in
  electricity pools.
\newblock {\em IEEE Transactions on Power Systems}, 18(2):452--459.

\bibitem[Gao and Li, 2011]{gao2011cardinality}
Gao, J. and Li, D. (2011).
\newblock Cardinality constrained linear-quadratic optimal control.
\newblock {\em IEEE Transactions on Automatic Control}, 56(8):1936--1941.

\bibitem[Goemans and Williamson, 1995]{goemans1995improved}
Goemans, M.~X. and Williamson, D.~P. (1995).
\newblock Improved approximation algorithms for maximum cut and satisfiability
  problems using semidefinite programming.
\newblock {\em Journal of the ACM (JACM)}, 42(6):1115--1145.

\bibitem[G{\"u}nl{\"u}k and Linderoth, 2010]{gunluk2010perspective}
G{\"u}nl{\"u}k, O. and Linderoth, J. (2010).
\newblock Perspective reformulations of mixed integer nonlinear programs with
  indicator variables.
\newblock {\em Mathematical Programming}, 124:183--205.

\bibitem[Han et~al., 2020]{HGA:2x2}
Han, S., G\'omez, A., and Atamt\"urk, A. (2020).
\newblock 2x2-convexifications for convex quadratic optimization with indicator
  variables.
\newblock {\em arXiv preprint arXiv:2004.07448}.

\bibitem[Hijazi et~al., 2012]{hijazi2012mixed}
Hijazi, H., Bonami, P., Cornu{\'e}jols, G., and Ouorou, A. (2012).
\newblock Mixed-integer nonlinear programs featuring “on/off” constraints.
\newblock {\em Computational Optimization and Applications}, 52:537--558.

\bibitem[Ho-Nguyen and K{\i}l{\i}n\c{c}-Karzan, 2017]{ho2017second}
Ho-Nguyen, N. and K{\i}l{\i}n\c{c}-Karzan, F. (2017).
\newblock A second-order cone based approach for solving the trust-region
  subproblem and its variants.
\newblock {\em SIAM Journal on Optimization}, 27(3):1485--1512.

\bibitem[Hochbaum, 2001]{hochbaum2001efficient}
Hochbaum, D.~S. (2001).
\newblock An efficient algorithm for image segmentation, {M}arkov random fields
  and related problems.
\newblock {\em Journal of the {ACM}}, 48:686--701.

\bibitem[Javanmard et~al., 2016]{javanmard2016phase}
Javanmard, A., Montanari, A., and Ricci-Tersenghi, F. (2016).
\newblock Phase transitions in semidefinite relaxations.
\newblock {\em Proceedings of the National Academy of Sciences},
  113(16):E2218--E2223.

\bibitem[Jeyakumar and Li, 2014]{jeyakumar2014trust}
Jeyakumar, V. and Li, G. (2014).
\newblock Trust-region problems with linear inequality constraints: exact {SDP}
  relaxation, global optimality and robust optimization.
\newblock {\em Mathematical Programming}, 147(1-2):171--206.

\bibitem[Pilanci et~al., 2015]{pilanci2015sparse}
Pilanci, P., Wainwright, M.~J., and El~Ghaoui, L. (2015).
\newblock Sparse learning via boolean relaxations.
\newblock {\em Mathematical Programming}, 151:63--87.

\bibitem[Shor, 1987]{shor1987quadratic}
Shor, N.~Z. (1987).
\newblock Quadratic optimization problems.
\newblock {\em Soviet Journal of Computer and Systems Sciences}, 25:1--11.

\bibitem[Wang and K{\i}l{\i}n{\c{c}}-Karzan, 2020]{wang2019generalized}
Wang, A.~L. and K{\i}l{\i}n{\c{c}}-Karzan, F. (2020).
\newblock The generalized trust region subproblem: solution complexity and
  convex hull results.
\newblock {\em Mathematical Programming}, pages 1--42.

\bibitem[Wang and K{\i}l{\i}n{\c{c}}-Karzan, 2021]{wang2019tightness}
Wang, A.~L. and K{\i}l{\i}n{\c{c}}-Karzan, F. (2021).
\newblock On the tightness of {SDP} relaxations of {QCQPs}.
\newblock {\em Mathematical Programming}, pages 1--41.

\bibitem[Wu et~al., 2017]{wu2017quadratic}
Wu, B., Sun, X., Li, D., and Zheng, X. (2017).
\newblock Quadratic convex reformulations for semicontinuous quadratic
  programming.
\newblock {\em {SIAM} Journal on Optimization}, 27:1531--1553.

\bibitem[Zheng et~al., 2014]{zheng2014improving}
Zheng, X., Sun, X., and Li, D. (2014).
\newblock Improving the performance of miqp solvers for quadratic programs with
  cardinality and minimum threshold constraints: A semidefinite program
  approach.
\newblock {\em INFORMS Journal on Computing}, 26(4):690--703.

\end{thebibliography}

\end{document}